\documentclass[a4paper,12pt,wlayout]{article}

\usepackage{array} 
\usepackage{amssymb}
\usepackage{amsmath} 
\usepackage{amsbsy} 
\usepackage{epsfig}
\usepackage{mathrsfs} 
\usepackage{graphics}
\usepackage{graphicx}
\usepackage[sort&compress,numbers]{natbib}
\usepackage{pstricks}
\usepackage{float}
\usepackage{amsthm}
\usepackage{slashed}
\usepackage{color}

\theoremstyle{plain}
\newtheorem{theorem}{\bf Theorem}[section]
\newtheorem{lemma}[theorem]{\bf Lemma}

\newtheorem{proposition}[theorem]{\bf Proposition}
\newtheorem{definition}[theorem]{\bf Definition}
\theoremstyle{remark}
\newtheorem{remark}[theorem]{\bf Remark}

\def\R{\mathrm{Re}}
\newcommand{\eps}{\varepsilon}
\newcommand{\lb}{\label}
\newcommand{\go}{\rightarrow}
\newcommand{\no}{\noindent}

\newcommand{\pt}{\partial_t}
\newcommand{\pxx}{\partial_{xx}}
\newcommand{\pxt}{\partial_{xt}}

\newcommand{\px}{\partial_x}

\newcommand{\ee}{\end{equation}}
\newcommand{\be}{\begin{equation}}
\newcommand{\bea}{\begin{eqnarray}}
\newcommand{\eea}{\end{eqnarray}}
\newcommand{\sbea}{\begin{subequations}\begin{eqnarray}}
\newcommand{\seea}{\end{eqnarray}\end{subequations}} 
\newcommand{\ees}{\end{equation*}}
\newcommand{\bes}{\begin{equation*}}
\newcommand{\beas}{\begin{eqnarray*}}
\newcommand{\eeas}{\end{eqnarray*}}

\newcommand{\rf}[1]{(\ref{#1})}

\newcommand{\myat}{\quad\mbox{at}\quad}
\newcommand{\const}{\mathrm{const}}
\newcommand{\mR}{\mathbb{R}\,}

\newcommand{\dv}{\mathrm{div}\,}

\usepackage{a4wide}
\begin{document}

\title{Weak solutions to lubrication equations in the presence of strong slippage}
\author{Georgy Kitavtsev\footnote{Max Planck Institute for Mathematics in the Sciences, Inselstr. 22, D--04103 Leipzig, Germany. {\tt E-mail}: Georgy.Kitavtsev@mis.mpg.de}, Philippe Lauren\c cot\footnote{Institut de Math\'ematiques de Toulouse, CNRS UMR~5219, Universit\'e de Toulouse, F--31062 Toulouse cedex 9, France. {\tt E-mail}: laurenco@math.univ-toulouse.fr},
and Barbara Niethammer \footnote{Mathematical Institute, University of Oxford, 24--29 St Giles', Oxford OX13LB, United Kingdom. {\tt E-mail}: niethammer@maths.ox.ac.uk}}
\date{\today}
\maketitle
\numberwithin{equation}{section}
\begin{abstract}
 The existence of global weak solutions is proved for one-dimensional
 lubrication models that describe the dewetting process of nanoscopic thin
 polymer films on hydrophobyzed substrates and take account of large slippage
 at the polymer-substrate interface. The convergence of these solutions
as either the Reynolds number or the capillarity goes to zero, 
as well as their limiting behaviour as the slip length goes to zero or infinity are investigated.
\end{abstract}

\section{Introduction.}
During the last thirty years lubrication theory was successfully applied to modeling of dewetting processes in micro and nanoscopic liquid films on a solid polymer substrates see e.g.~\citet{ODB97,MWW06} and references therein. The influence of intermolecular interactions that are typically due to the competition between the long-range attractive van der Waals and short-range Born repulsive intermolecular forces play an important role in such processes, see~\citet{gennes85,ODB97}.

Besides intermolecular forces and surface tension at the free surface of the film it has been shown by~\citet{FMWJ07} that the dewetting of polymer films on hydrophobic substrates also involves such boundary effect as slippage on a solid substrate. The measure of slip is a so-called slip length, which is defined as an extrapolated distance relative to the wall where the tangential velocity component of the liquid vanishes. Recently, it has been shown experimentally and theoretically that the early stages of the dewetting process and the evolving morphology depend markedly on the magnitude of the effective slip length, which can be of the size of the height of the liquid film or even larger for nanoscale systems, see e.g.~\citet{FJMWW05,MW05}. Recently in~\citet{MWW06,KSK04} closed-form one-dimensional lubrication equations over a wide range of slip lengths were derived from the underlying equations for conservation of mass and momentum, together with boundary conditions for the tangential and normal stresses, as well as the kinematic condition at the free boundary, impermeability and Navier-slip condition at the liquid-solid interface. Asymptotic arguments, based on the magnitude of the slip length show that within a lubrication scaling there are two \emph{distinguished limits}, see~\citet{MWW06}.

These are the well-known \emph{weak-slip} model 
\be
\pt h= - \px \Big((h^3+b\,h^2)\px \left(\sigma \pxx h-\Pi(h)\right)\Big) 
\label{WSM}
\ee
with $b$ denoting the slip length parameter, and the \emph{strong-slip} model 
\begin{subequations}
\label{SSM}
\begin{align}
\R\,(\pt(hu)+\px(hu^2))&=4\px(\nu(h)\px u)+h\px(\sigma\pxx h-\Pi(h))-\frac{u}{\beta}
\label{SSM1}\\
\pt h &= -\,\px\left(h u\right),
\label{SSM2}
\end{align}
\end{subequations}
respectively. Here, $u(x,t)$ and $h(x,t)$ denote the average velocity in the lateral direction and the height profile for the free surface, respectively. The positive slip length parameters $b$ and $\beta$ are related by orders of magnitude via $b\sim\eta^2\beta$, where the (small) parameter $\eta$, $0<\eta\ll 1$, refers to the vertical to horizontal scale separation of the thin film. The high order of the lubrication equations \rf{WSM} and \rf{SSM1}--\rf{SSM2} is a result of the contribution from surface tension at the free boundary, reflected by the linearized curvature term $\sigma \pxx h$ with parameter $\sigma\ge 0$.  A further contribution to the pressure is denoted by $\Pi(h)$ and represents that of the intermolecular forces, namely long-range attractive van der Waals and short-range Born repulsive intermolecular forces. A commonly used expression for it is given by
\be
\label{PF} \Pi(h)=\frac{1}{h^3}-\frac{\alpha}{h^4}\ \ \text{with}\ \ \alpha>0. 
\ee
The terms $\R\,(\pt(hu)+\px(hu^2))$, with $\R\ge 0$ denoting the Reynolds number, and $4\px(\nu(h)\px u)$ in \rf{SSM1}--\rf{SSM2} represent inertial and Trouton viscosity terms, respectively. We assume below linear dependence of the viscosity coefficient on the height, namely
\bes
\nu(h):=\nu h,\ \nu>0.
\lb{vs}
\ees
 
\noindent For these choice of the constitutive laws, we first investigate the existence of global weak solutions to \rf{SSM1}--\rf{SSM2} on the interval $(0,1)$, supplemented with the boundary conditions
\begin{eqnarray}
u &=& 0 \myat x=0,\,1\,,\label{SSM flux BC} \\
\px h &=& 0 \myat x=0,\,1\,,\label{SSM height BC}
\end{eqnarray}
and initial data $(h_0,u_0)$ with positive first component $h_0$. Observe that, in particular, the boundary conditions \rf{SSM flux BC} for $u$ guarantee the conservation of mass 
\be
\int_{0}^{1}h(x,t)\,dx = h_c := \int_0^1 h_0(x)\, dx\,, \quad  t\ge 0,
\lb{conservation of mass}
\ee
where $h_c$ is the average of the height profile.
 
\noindent  We next investigate the behaviour of these weak solutions as some of the parameters in the model goes to zero or to infinity. More precisely, let us point out that the strong-slip model \rf{SSM1}--\rf{SSM2} includes as limiting cases three further lubrication models. One is obtained from the strong-slip model in the limit $\beta\rightarrow\infty$ and describes the  dynamics of \emph{suspended free films}, see e.g.~\citet{BG99}:
\begin{subequations}
\label{FFM}
\begin{align}
\R\,(\pt(hu)+\px(hu^2))&=4\px(\nu(h) \px u)+h\px(\sigma\pxx h-\Pi(h)),
\label{FFM1}\\
\pt h &= -\,\px\left(h u\right).
\label{FFM2}
\end{align}
\end{subequations}
 In the second one, we neglect the inertial terms which corresponds to a vanishing Reynolds number. For the third limiting case which is derived in~\citet{MWW06} the slip length parameter $\beta_I$ is of order of magnitude lying in between those that lead to the weak-slip model \rf{WSM} and the strong-slip model \rf{SSM1}--\rf{SSM2}, i.e. $b\ll \beta_I\ll\beta$. The corresponding \emph{intermediate-slip} model is given by 
\be
\pt h= - \px \Big(h^2\px \left(\sigma \pxx h-\Pi(h)\right)\Big) \,.
\label{ISM}
\ee
It can be obtained by rescaling time in \rf{WSM} by $b$ and letting $b\rightarrow\infty$ or by rescaling time and the horizontal velocity by $\beta$ in \rf{SSM1}--\rf{SSM2} and taking the limit $\beta\rightarrow 0$. We consider here the second limit and change variables in \rf{SSM1}--\rf{SSM2}  as follows:
\be
x:= \bar{x}\,, \ \ t:=\bar{t}\beta\,, \ \ h:=\bar{h}\,, \ \ \text{and}\ \ u:=\frac{\bar{u}}{\beta},
\lb{CoV}
\ee
where $\bar{x}$, $\bar{t}$, $\bar{h}$, and $\bar{u}$ denote the old variables. We show that, as $\beta\go 0$, the rescaled solutions to \rf{SSM1}--\rf{SSM2} converge to a solution to
\begin{subequations}
\label{ISMl}
\begin{align}
u&=h\px(\sigma\pxx h-\Pi(h)),
\label{ISM1}\\
\pt h &= -\,\px\left(h u\right),
\label{ISM2}
\end{align}
\end{subequations}
satisfying boundary conditions \rf{SSM flux BC}--\rf{SSM height BC} that is, a solution to the intermediate-slip equation \rf{ISM}.

Our proofs follow a strategy that has also been employed in other systems
arising in fluid dynamics and phase transition theories.
In fact,  \citet{BDL03} and \citet{Br09} considered the following 
Korteweg system 
\begin{subequations}
\label{GKM}
\begin{align}
\pt(hu)+\dv(hu\otimes u)&=\dv(\mu(h)(\nabla u+(\nabla u)^T)+\nabla(\lambda(h)\dv u)
\label{GKM1}\\
&+\sigma h \nabla\Delta h-\nabla P(h) \nonumber\\
\pt h &= -\,\dv\left(h u\right).
\label{GKM2}
\end{align}
\end{subequations}  
In the case  
\bes
\sigma=0,\,\mu(h)=\mu=\const,\,\lambda(h)=\lambda=\const ,
\ees
the system\rf{GKM1}--\rf{GKM2} gives the
 compressible Navier-Stokes equations. For the latter, considered on a bounded domain with $P(h)$ having a globally
Lipschitz continuous derivative,  existence of a unique strong solution was shown by~\citet{So76} in dimensions $d\ge 2$. 
Recently,~\citet{BDL03} showed existence of a weak solution to \rf{GKM1}--\rf{GKM2}  on  bounded domain with 
\be
\sigma\ge 0,\,\mu(h)=\nu h,\,\lambda(h)=0
\lb{M6}
\ee
and similar requirements on $P(h)$ for dimensions $d\ge 2$. In the proof the authors use a so
called BD entropy equality introduced by them originally in~\citet{BD03} for
the case $\sigma=0$. In the latter case \rf{GKM1}--\rf{GKM2} coincides with
viscous shallow-water equations. In a recent review by~\citet{Br09} on those
equations it was shown that the BD entropy equality  holds for \rf{GKM1}--\rf{GKM2}  in
dimensions $d\ge 2$ if a special relation between
viscosity coefficients is satisfied, namely
\be
\lambda(h)=2(\mu'(h)h-\mu(h))
\lb{lmR}
\ee
that is satisfied for \rf{M6}.

We also notice that the strong-slip lubrication equation \rf{SSM1}--\rf{SSM2} has a similar form to
the one dimensional case of \rf{GKM1}--\rf{GKM2}. Indeed, in one space dimension
\rf{GKM1}--\rf{GKM2} reduces to
\begin{subequations}
\label{KM}
\begin{align}
\pt(hu)+\px(hu^2)&=\px(\nu(h)\px u)+\sigma h\px^3 h-\px P(h)
\label{KM1}\\
\pt h &= -\,\px\left(h u\right).
\label{KM2}
\end{align}
\end{subequations}
Recently, in~\citet{MV08} an analogue of the BD entropy identity was shown for \rf{KM1}--\rf{KM2} with  any $C^2$ smooth viscosity function $\nu(h)$ and $\sigma=0$. Furthermore,  they proved existence of strong solutions to \rf{KM1}--\rf{KM2}  on the whole real line $\mR$ with $\sigma=0,\,\nu(h)=\nu h^{k}\ \ \text{and}\ \ P(h)=h^\gamma$, where $k\le 1/2$ and $\gamma>1$.

\noindent The main difference in \rf{SSM1}--\rf{SSM2} compared to \rf{KM}  is the special form of the intermolecular pressure \rf{PF} which is singular when the height $h$ vanishes and complicates the analysis. Nevertheless, it also provides a very useful positive lower bound for $h$. Equations \rf{SSM1}--\rf{SSM2} also include an additional slip term which is easily handled for the existence theory but plays a crucial role in the limit $\beta\go 0$ (leading to equation \rf{ISM}).

In this paper we start in section~\ref{S.apriori} with recalling the energy identity associated to \rf{SSM1}--\rf{SSM2} and derive  a version of a BD entropy identity.  These two estimates provide the main ingredients of our existence results. In section~\ref{S.existence} we follow the strategy of~\citet{BD06} by setting up a higher-order regularized equation for which analogous energy and entropy inequalities are satisfied. The a priori estimates allow us to pass to the limit in the regularization parameter to obtain global weak solutions to \rf{SSM1}--\rf{SSM2}. In section~\ref{S.limit1} we establish existence of global weak solutions for the cases
$\R=0, \beta=\infty$ and $\sigma=0$ by passing to the respective limits in \rf{SSM1}--\rf{SSM2}. While the first
two cases are straightforward, we need to refine the estimates from section~\ref{S.apriori} in order to enable
us to deal with the limit $\sigma \to 0$. Finally, in section \ref{S.limit2} we establish that solutions of
\rf{SSM1}--\rf{SSM2} converge, after the rescaling \rf{CoV},  in the limit $\beta \to 0$ to a solution of the intermediate-slip
model \rf{ISM}. 

\section{A priori estimates.}\label{S.apriori}

In this section we state two relations that are satisfied by classical solutions of \rf{SSM1}--\rf{SSM2}. As it was stated in the introduction both have their counterparts for viscous shallow-water equation \rf{GKM1}--\rf{GKM2} and lubrication equations \rf{WSM} and \rf{ISM}. For the latter they were initially derived by~\citet{BF90}. The energy equality for shallow-water equations is known already for several decades, whereas the entropy equality was suggested recently by~\citet{BD03}. For the strong-slip model \rf{SSM1}--\rf{SSM2} the derivation of the energy equality is again standard, see~\citet{KW09}. As a new result we derive here the entropy equality  for \rf{SSM1}--\rf{SSM2} following the approach of~\citet{BD03}. At the end of this section we use the energy and entropy equalities to derive a  priori estimates in the case $\sigma>0$ on classical solutions to \rf{SSM1}--\rf{SSM2}  having positive $h$ and satisfying \rf{conservation of mass}, \rf{SSM flux BC}--\rf{SSM height BC}. For consistency we start with the energy equality.
\begin{lemma}[Energy equality]
For classical solutions of the system \rf{SSM1}--\rf{SSM2} with boundary conditions \rf{SSM flux BC}--\rf{SSM height BC} the following equality holds
\be
\frac{d}{dt}\int_0^1\left[\R\,h\frac{u^2}{2}+U(h)+\sigma\frac{|\px h|^2}{2}\right]\,dx=-4\int_0^1\nu h|\px u|^2\, dx-\int_0^1\frac{u^2}{\beta}\,dx,
\lb{EE}
\ee
where the potential function $U$ is the indefinite integral of $\Pi$ defined by 
\be
U(h)=-\frac{1}{2\,h^2}+\frac{\alpha}{3\,h^3}\,, \quad h>0\,.
\label{U}
\ee
\lb{L1}
\end{lemma}
\no\textbf{Proof:} One way to show \rf{EE} is to use the fact that 
\be
E(u,h):=\int_0^1\left[\R\, h\frac{u^2}{2}+U(h)+\sigma\frac{|\px h|^2}{2}\right]\,dx
\lb{LF}
\ee
is a Lyapunov functional for the system \rf{SSM1}--\rf{SSM2}, see~\citet{KW09}. Here we give another standard derivation of \rf{EE}. Multiplying \rf{SSM1} by $u$ and integrating in $x$ we obtain
\beas
&&\R\,\int_0^1\left[\pt(hu)u+\px(hu^2)u\right]\,dx-4\int_0^1\px(\nu\,h\px u)u\,dx\\
&&-\int_0^1hu\px(\sigma\pxx h-\Pi(h))\,dx+\int_0^1\frac{u^2}{\beta}\,dx=0.
\eeas
Using now several integrations by parts, \rf{SSM2} and boundary conditions \rf{SSM flux BC}--\rf{SSM height BC} we obtain
\beas
0&=&\R\,\int_0^1\left[ \pt\left( h\frac{u^2}{2}\right)+\pt h\frac{u^2}{2}-h u^2\px u \right]\,dx\\
&+&4\int_0^1\nu\,h(\px u)^2\,dx-\int_0^1\pt h(\sigma\pxx h-\Pi(h))\,dx+\int_0^1\frac{u^2}{\beta}\,dx\\
&=&\R\,\int_0^1\left[ \pt\left( h\frac{u^2}{2}\right)-\px (h u)\frac{u^2}{2}-hu\px\left(\frac{u^2}{2}\right) \right]\,dx\\
&+&4\int_0^1\nu\,h(\px u)^2\,dx+\int_0^1\pt \left(\sigma\frac{|\px h|^2}{2}+U(h)\right)\,dx+\int_0^1\frac{u^2}{\beta}\,dx\\
&=&\R\,\int_0^1\pt\left( h\frac{u^2}{2}\right)\,dx+4\int_0^1\nu\,h(\px u)^2\,dx+\int_0^1\pt \left(\sigma\frac{|\px h|^2}{2}+U(h)\right)\,dx
+\int_0^1\frac{u^2}{\beta}\,dx.
\eeas
From the last expression \rf{EE} follows.$\hfill\blacksquare\\$
\begin{lemma}
For classical solutions of the system \rf{SSM1}--\rf{SSM2} with boundary conditions \rf{SSM flux BC}--\rf{SSM height BC} the following equality holds
\be
\frac{1}{2}\frac{d}{dt}\int_0^1\frac{|\px h|^2}{h}\,dx=\int_0^1\px\left(\frac{\px h}{h}\right)\px u\,h\,dx.
\lb{P1}
\ee
\lb{Pr1}
\end{lemma}
\no\textbf{Proof:} The statement is verified again using several integrations by parts, equation \rf{SSM1} and boundary conditions \rf{SSM flux BC}--\rf{SSM height BC} as follows.
\beas
&&\frac{1}{2}\frac{d}{dt}\int_0^1\frac{|\px h|^2}{h}\,dx=\int_0^1\frac{\px h\pxt h}{h}\,dx-\frac{1}{2}\int_0^1\frac{\pt h |\px h|^2}{h^2}\,dx\\
&=&-\int_0^1\frac{\px h}{h}\pxx(hu)\,dx+\frac{1}{2}\int_0^1\px (hu)\left(\frac{\px h}{h}\right)^2\,dx\\
&=&\int_0^1\px\left(\frac{\px h}{h}\right)\px(hu)\,dx-\int_0^1u \px h\px\left(\frac{\px h}{h}\right)\,dx\\
&=&\int_0^1\px\left(\frac{\px h}{h}\right)\px u\,h\,dx.
\eeas$\hfill\blacksquare\\$
Let us define a so-called entropy function by $\varphi(h):=4\nu\log h$ for $h>0$. Then we have the following lemma.

\begin{lemma}[Entropy equality]
For classical solutions of the system \rf{SSM1}--\rf{SSM2} with boundary conditions \rf{SSM flux BC}--\rf{SSM height BC} the following equality holds.
\bea
&&\frac{d}{dt}\int_0^1\left[ \frac{1}{2}h(\R\, u+\px \varphi(h))^2  -\frac{1}{\beta}\varphi(h) +\R \left( \sigma\frac{|\px h^2}{2}+U(h) \right) \right]\,dx\phantom{0000000}\nonumber\\
&&=-\R\,\int_0^1\frac{u^2}{\beta}\,dx-4\sigma\nu\int_0^1|\pxx h|^2\,dx-4\nu\int_0^1\Pi'(h)|\px h|^2\,dx.
\lb{EnE}
\eea
\lb{L2}
\end{lemma}
\no\textbf{Proof:} 
Let us multiply equation \rf{SSM1} by $\px \varphi(h)$, integrate with respect to $x$ and use an integration by parts, \rf{SSM2} and \rf{SSM height BC}:
\bea
&&4\nu \R\,\int_0^1(\pt u+u\px u)\px h\,dx=-16\nu^2\int_0^1\px\left(\frac{\px h}{h}\right)\px u\,h\,dx\nonumber\\
&&-4\nu\int_0^1\frac{u\px h}{\beta h}\,dx-4\sigma\nu\int_0^1|\pxx h|^2\,dx-4\nu\int_0^1\Pi'(h)|\px h|^2\,dx.
\lb{P2}
\eea
On one hand, a further integration by  parts in the left-hand side of \rf{P2},  equation \rf{SSM2}, and the energy equality \rf{EE} give
\beas
&&4\nu \R\,\int_0^1(\pt u+u\px u)\px h\,dx\\
&=&4\nu \R\,\left( \frac{d}{dt}\int_0^1u\px h\,dx-\int_0^1u\pxt h\,dx+\int_0^1u\px u\px h\,dx \right)\\
&=&4\nu \R\,\left( \frac{d}{dt}\int_0^1u\px h\,dx-\int_0^1\px u\px(uh)\,dx+\int_0^1u\px u\px h\,dx \right)\\
&=&4\nu \R\,\left( \frac{d}{dt}\int_0^1u\px h\,dx-\int_0^1h(\px u)^2\,dx \right)\\
&=&\R\,\left(\frac{d}{dt} \int_0^1 \left[ 4\nu u\px h+\R\, h\frac{u^2}{2}+U(h)+\sigma\frac{|\px h|^2}{2} \right]\,dx +\int_0^1\frac{u^2}{\beta}\,dx\right).
\eeas
On the other hand, using Lemma~\ref{Pr1} and the definition of $\varphi$, one can write the first term on the right-hand side of \rf{P2} as
\bes
-16\nu^2\int_0^1\px\left(\frac{\px h}{h}\right)\px u h\,dx=-\frac{1}{2}\frac{d}{dt}\int_0^1h\left|\px \varphi(h)\right|^2\,dx,
\ees
while it follows from \rf{SSM2} and \rf{SSM flux BC} that the second term on the right-hand side of \rf{P2} can be transformed into
\bes
-4\nu\int_0^1\frac{u\px h}{\beta h}\,dx=-4\nu\int_0^1\frac{\px(uh)}{\beta h}\,dx+4\nu\int_0^1\frac{\px u}{\beta}\,dx =\frac{1}{\beta}\frac{d}{dt}\int_0^1\varphi(h)\,dx.
\ees
Substituting finally the last three identities into \rf{P2} one obtains
\beas
&&\frac{d}{dt} \int_0^1 \left[ \frac{\R\,^2}{2}hu^2+\frac{1}{2}h\left|\px \varphi(h)\right|^2 +4\nu \R\,u\px h-\frac{1}{\beta}\varphi(h)+\R\,\left( \sigma\frac{|\px h|^2}{2}+U(h) \right) \right] \,dx\\
&&=-\R\,\int_0^1\frac{u^2}{\beta}\,dx-4\sigma\nu\int_0^1|\pxx h|^2\,dx-4\nu\int_0^1\Pi'(h)|\px h|^2\,dx.
\eeas
Using the definition of $\varphi$ the last expression can be easily transformed into \rf{EnE}.$\hfill\blacksquare\\$

\begin{lemma}
For smooth functions $(h,u)$ with a positive first component $h$, we have
\be
\frac{1}{4}\int_0^1h|\px\varphi(h)|^2\,dx\le \frac{1}{2}\int_0^1h\left(\R\, u+\px\varphi(h)\right)^2\,dx+\R\,E(u,h)+\frac{\R}{6\alpha^2},
\lb{CoEs}
\ee
where $E(u,h)$ is defined in \rf{LF}.
\lb{C1}
\end{lemma}
\no\textbf{Proof:}
Using the elementary inequality 
\bes
(y+z)^2\ge\frac{y^2}{2}-z^2,
\ees
the fact that 
\be
U(h)\ge -1/(6\alpha^2)
\lb{UBd}
\ee
for all $h>0$, and the definition \rf{LF} of $E$ one obtains
\beas
&&\int_0^1h\left(\R\, u+\px\varphi(h)\right)^2\,dx\ge\frac{1}{2}\int_0^1h|\px\varphi(h)|^2\,dx-\R^2\,\int_0^1h u^2\,dx\\
&\ge&\frac{1}{2}\int_0^1h|\px\varphi(h)|^2\,dx-2\R\left[E(u,h)-\int_0^1\left(U(h)+\sigma\frac{|\px h|^2}{2}\right)\,dx\right]\\
&\ge&\frac{1}{2}\int_0^1h|\px\varphi(h)|^2\,dx-2\R\, E(u,h)-\frac{\R}{3\alpha^2},
\eeas
from which the statement of the lemma follows.$\hfill\blacksquare\\$

\bigskip

\no We now deduce several bounds from Lemmas~\ref{L1} and~\ref{L2} for classical solutions $(h,u)$ to the system \rf{SSM1}--\rf{SSM2} with boundary conditions \rf{SSM flux BC}--\rf{SSM height BC} and a positive first component. To this end, we assume that the initial data $h_0$ and $u_0$ for \rf{SSM1}--\rf{SSM2} satisfy
\bea
h_0\in H^1(0,1),\ \ h_0(x)>0\ \ \text{for all}\ \ x\in [0,1]\ \ \text{and}\ \ u_0\in L^2(0,1).
\lb{IAs}
\eea
As already mentioned, a first consequence of \rf{SSM2} and \rf{SSM flux BC} is the conservation of mass \rf{conservation of mass}.
We next combine \rf{conservation of mass}, the energy and entropy identities to obtain the following bounds when $\sigma>0$.

\begin{proposition}[a priori estimates]
For fixed positive $\sigma,\R,\beta,T$  and initial data satisfying \rf{IAs} there exists a positive constant $C_0>1$ depending only on $T$, $\alpha$, $\nu$, $\R$, $\sigma$, $E(u_0,h_0)$, and $S_\beta(u_0,h_0)$ defined in \rf{Sb} below such that the following terms are bounded by $C_0$ in respective norms
\bea
\sqrt{h},\,\px \sqrt{h},\,h^{-3/2},\, \px h,\,\sqrt{\R}\sqrt{h}u&\in& L^\infty(0,T;L^2(0,1)),\lb{AE0s}\\
\px(h^{-3/2}),\,\pxx h,\,\sqrt{h}\px u,\,\frac{u}{\sqrt{\beta}}&\in& L^2((0,1)\times (0,T)),
\lb{AEs}
\eea
and
\be
C_0^{-1} \le h(x,t)\le C_0
\lb{BEs}
\ee
for all $x\in (0,1)$ and $t\in (0,T)$.  The constant $C_0$ tends to $\infty$ as $\sigma\go 0$. 
\lb{Pr2}
\end{proposition}
\no\textbf{Proof:}  
Integrating the energy equality \rf{EE} with respect to time over $(0,\,t)$,  $t\in (0,\,T)$, and using the inequality
\bes
U(h)\ge \frac{\alpha}{6h^3}-\frac{2}{3\alpha^2}
\ees
ensure that
\beas
\|\sqrt{\R} \sqrt{h(t)} u(t)\|_2^2 + \frac{\alpha}{3} \|h^{-3/2}(t)\|_2^2 + \|\sqrt{\sigma} \px h(t)\|_2^2 &\le & C_1:= 2 E(u_0,h_0) + \frac{4}{3\alpha^2}, \\
\int_0^t \left[ 4\nu \|\sqrt{h(s)} \px u(s)\|_2^2 + \left\| \frac{u(s)}{\sqrt{\beta}} \right\|_2^2\right]\,ds &\le & C_1.
\eeas
Therefore, the norms of $\sqrt{\sigma}\px h$, $h^{-3/2}$, and $\sqrt{\R}\sqrt{h}u$ in $L^\infty(0,T;L^2(0,1))$ are bounded by a constant that depends only on $C_1$ while the norms of $\sqrt{h}\px u$ and $u/\sqrt{\beta}$ in $L^2((0,1)\times (0,T))$ are bounded by a constant that depends only on $C_1$ and $\nu$.Using the bound on $\px h$ one obtains
$$
|h(x,t)-h(y,t)|=\left|\int_x^y\px h(z,t)\,dz\right|\le |x-y|^{1/2}\ \|\px h(t)\|_2 \le \frac{C_1}{\sqrt{\sigma}}|x-y|^{1/2}
$$
for all $(x,y)\in (0,1)\times (0,1)$ and $t\in (0,T)$. Integrating the above inequality with respect to $y\in (0,1)$ and using  \rf{conservation of mass} readily give the upper bound in\rf{BEs}. To establish the lower bound for $h$ in \rf{BEs}, we combine the $L^\infty(0,T;L^2(0,1))$-estimates on $h^{-3/2}$ and $\sqrt{\sigma}\px h$ just established to obtain a bound on the norm of $1/\sqrt{h}$ in $L^\infty(0,T;W^{1,1}(0,1))$ since
$$
\int_0^1 \left| \px\left( h^{-1/2} \right) \right| \,dx = \frac{1}{2} \int_0^1 \frac{ \left| \px h \right|}{h^{3/2}}\,dx \le \frac{1}{2\sqrt{\sigma}}\ \|\sqrt{\sigma}\px h\|_2\ \|h^{-3/2}\|_2.
$$
Owing to the continuous embedding of $W^{1,1}(0,1)$ in $L^\infty(0,1)$, the positive lower bound in \rf{BEs} readily follows.

Next, let us introduce the functional $S_\beta(u,h)$ defined by
\be
S_\beta(u,h):=\int_0^1\Big[\frac{1}{2}h(\R\, u+\px \varphi(h))^2+\frac{1}{\beta}(4\nu h-\varphi(h)) +\R\,\sigma\frac{|\px h|^2}{2}+\R\,U(h)\Big]\,dx.
\lb{Sb}
\ee
It follows from the mass conservation \rf{conservation of mass} and the entropy equality \rf{EnE} that
\be
S_\beta(u(t),h(t))+\int_0^t\int_0^1\left(\frac{\R}{\beta}u^2+4\sigma\nu|\pxx h|^2+4\nu\Pi'(h)|\px h|^2\right)\,dx\,ds=S_\beta(u_0,h_0)
\lb{SEq}
\ee
for any $t\in(0,T)$. Since $z-\log(z)\ge 1$ for all $z>0$,  it follows from the energy equality \rf{EE}, \rf{UBd}, and Lemma~\ref{C1} that
\bea
S_\beta(u,h)&\ge & \frac{1}{2}\int_0^1h(\R\, u+\px\varphi(h))^2\,dx+\frac{1}{\beta}\int_0^1(4\nu h-\varphi(h))\,dx-\frac{\R}{6\alpha^2}\nonumber\\
&\ge &\frac{1}{2}\int_0^1h(\R\, u+\px\varphi(h))^2\,dx-\frac{\R}{6\alpha^2} \lb{contrex}\\
&\ge & \frac{1}{4} \int_0^1 h |\px\varphi(h)|^2\,dx - \R\, E(u,h) - \frac{\R}{3\alpha^2} \nonumber\\
&\ge & \frac{1}{4} \int_0^1 h |\px\varphi(h)|^2\,dx - \R\, E(u_0,h_0) - \frac{\R}{3\alpha^2}. \nonumber
\eea
Combining this with the previous estimate on $\sqrt{\sigma}\px h$, \rf{SEq}, and the fact that 
\bes
\Pi'(h)\ge \frac{2\alpha}{h^5}-\left(\frac{6}{5}\right)^5\frac{1}{2\alpha^4}
\ees
implies that
\beas
& & \frac{1}{4} \int_0^1 h(t) |\px\varphi(h(t))|^2\,dx + \int_0^t \int_0^1 \left( 4\sigma\nu|\pxx h(s)|^2+\frac{8\nu\alpha|\px h(s)|^2}{h^5(s)}\right)\,dx\,ds\\
& \le & C_2 + \left(\frac{6}{5}\right)^5\frac{2\nu}{\alpha^4} \int_0^t \int_0^1 |\px h(s)|^2\,dx\,ds \le C_2 + \left(\frac{6}{5}\right)^5\frac{2\nu}{\alpha^4} \frac{C_1^2 T}{\sigma}
\eeas
with $C_2:=S_\beta(u_0,h_0)+ \R\,E(u_0,h_0) +(\R/3\alpha^2)$. Since $h|\px\varphi(h)|^2=4 |\px\sqrt{h}|^2$, 
this completes the proof of \rf{AE0s}--\rf{AEs}.  $\hfill\blacksquare\\$

\begin{remark}
As in Theorem~1 of~\citet{BGW01} that shows existence of uniform lower bound for solutions of~\rf{ISM} many estimates in the proof of Proposition~\ref{Pr2} heavily rely on the positivity of the parameter $\alpha$ in $\Pi$ and become unbounded as  $\alpha\go 0$. The estimates \rf{AE0s}, \rf{AEs}, and \rf{BEs} are no longer valid in that case. The same remark holds for the case
$\sigma=0$ since the estimates depend on the bound  for $\sqrt{\sigma}\px h$. On the contrary the proof holds without changes for the cases $\R=0$ or $\beta=\infty$ except for the fact that then one looses estimates on $\R\sqrt{h} u$ and $u/\sqrt{\beta}$ in \rf{AE0s}--\rf{AEs}, respectively.
\lb{R2} 
\end{remark} 

\section{Existence of weak solutions.}
\label{S.existence}

We first define a weak formulation of the problem \rf{SSM1}--\rf{SSM2} with boundary conditions \rf{SSM flux BC}--\rf{SSM height BC}. Consider two functions $h_0$ and $u_0$ satisfying \rf{IAs}. 

\begin{definition}\lb{defws}
A pair $(h,u)$ is a global weak solution to \rf{SSM1}--\rf{SSM2} with boundary conditions \rf{SSM flux BC}--\rf{SSM height BC} and initial conditions $(h_0,u_0)$  if $h$ and $u$ enjoy the regularity properties stated in \rf{AE0s}--\rf{AEs}--\rf{BEs} and the following holds
\bea
\lb{WF1}
&&\int_0^\infty\int_0^1h\pt\psi\,dxdt + \int_0^1 h_0 \psi(\cdot,0)\,dx=-\int_0^\infty\int_0^1hu\px\psi\,dxdt,\\[2ex]
&&\R\,\int_0^\infty\int_0^1hu\pt\phi\,dxdt+\R\,\int_0^1h_0 u_0 \phi(\cdot,0)\,dx+\R\,\int_0^\infty\int_0^1hu^2\px\phi \, dx\, dt \nonumber\\
&-&4\nu\int_0^\infty\int_0^1h\px u\px\phi\,dx\,dt - \sigma\int_0^\infty\int_0^1\px h\pxx h \phi \,dx\,dt - \sigma \int_0^\infty\int_0^1h\pxx h\px\phi\,dx\,dt \nonumber\\
&+&\int_0^\infty\int_0^1\Pi_1(h)\px\phi\,dxdt-\frac{1}{\beta}\int_0^\infty\int_0^1u \phi\,dxdt=0
\lb{WF2}
\eea
for all $\psi\in C_0^\infty([0,1]\times [0,\infty))$ and $\phi\in C_0^\infty((0,1)\times [0,\infty))$, where
\be
\Pi_1(h):=-\int_h^\infty\tau\Pi'(\tau)\,d\tau.
\lb{Pi1}
\ee
\end{definition}

To show the existence of weak solutions to \rf{SSM1}--\rf{SSM2}--\rf{SSM flux BC}--\rf{SSM height BC}  we construct approximating systems similar to those suggested by~\citet{BD06} for the Korteweg and viscous shallow-water equations \rf{GKM1}--\rf{GKM2}. Although the pressure term $\Pi(h)$ does not need a regularization as in the case of~\citet{BD06}, one still needs to regularize the function $h$ sufficiently in order to control additional higher order terms arising in the entropy equality. The approximating systems we take are given by
\begin{subequations}
\label{ASSM}
\begin{align}
\R\,(\pt(h_{\eps}u_{\eps})+\px(h_{\eps}u_{\eps}^2))&=4\px(\nu\,h_{\eps}\px u_{\eps})+h_{\eps}\px(\sigma\pxx h_{\eps} - \Pi(h_{\eps})) - \frac{u_{\eps}}{\beta}\label{ASSM1}\\
& \quad +\eps h_{\eps}\partial_x^7 h_{\eps}-\eps^2\partial_x^4 u_{\eps},\phantom{0000}\nonumber\\
\pt h_{\eps} &= -\,\px\left(h_{\eps} u_{\eps}\right),
\label{ASSM2}
\end{align}
\end{subequations}
where $\eps>0$ is a small parameter. Consider \rf{ASSM1}--\rf{ASSM2} with boundary conditions
\be
u_{\eps}=\pxx u_{\eps}=\px h_{\eps}=\partial_x^3 h_{\eps}=\partial_x^5 h_{\eps}=0, \quad (x,t)\in \{0,1\}\times (0,T),
\label{ASSM BC}
\ee
and initial conditions
\be
u_{\eps}(x,0)=u_{\eps,0}(x)\ \ \text{and}\ \ h_{\eps}(x,0)=h_{\eps,0}(x)>0, \quad x\in (0,1),
\label{ASSM IC}
\ee
where $u_{\eps,0}$ and $h_{\eps,0}$ are smooth functions such that
\be
u_{\eps,0}\go u_{0}\ \text{in}\ L^2(0,1), \ h_{\eps,0}\go h_{0}\ \text{in}\ H^1(0,1), \ \text{and}\ \ \eps h_{\eps,0}\go 0\ \text{in}\ H^3(0,1)\ \text{as}\ \eps\go 0.
\lb{AP}
\ee
Proceeding as in Lemmas~\ref{L1} and~\ref{L2} one derives formulas analogous to \rf{EE} and \rf{EnE} for the system \rf{ASSM1}--\rf{ASSM2} with boundary and initial conditions \rf{ASSM BC}--\rf{ASSM IC}. Namely, introducing the energy
$$
E_\eps(u_\eps,h_\eps) := \int_0^1 \left[ \R\, h_{\eps}\frac{u_{\eps}^2}{2}+U(h_{\eps}) +\sigma\frac{|\px h_{\eps}|^2}{2}+\frac{\eps}{2}|\partial_x^3 h_{\eps}|^2 \right] \,dx,
$$
and the entropy
\beas
S_\eps(u_\eps,h_\eps) &:=& \int_0^1 \left[ \frac{1}{2} h_\eps (\R\, u_{\eps}+\px\varphi(h_{\eps}))^2 -\frac{1}{\beta} \varphi(h_{\eps}) \right]\,dx \\
&+& \R \int_0^1 \left[ \sigma\frac{|\px h_{\eps}|^2}{2}+U(h_{\eps}) + \frac{\eps}{2}|\partial_x^3 h_{\eps}|^2 \right] \,dx,
\eeas
the corresponding energy and entropy equalities read 
\be
\frac{d}{dt} E_\eps(u_\eps,h_\eps) = -4\int_0^1\nu h_{\eps}|\px u_{\eps}|^2\, dx - \int_0^1\frac{u_{\eps}^2}{\beta}\,dx -\eps^2\int_0^1|\pxx u_{\eps}|^2\,dx \lb{AEE}
\ee
and 
\bea
\frac{d}{dt} S_\eps(u_\eps,h_\eps) &=&-\R\,\int_0^1\frac{u_{\eps}^2}{\beta}\,dx-4\sigma\nu\int_0^1|\pxx h_{\eps}|^2\,dx-4\nu\int_0^1\Pi'(h_{\eps})|\px h_{\eps}|^2\,dx\nonumber\\
&-&\int_0^1\left[ \R\,\eps^2 |\pxx u_{\eps}|^2+4\nu\eps |\partial_x^4 h_{\eps}|^2 + 4\nu\eps^2 \pxx u_{\eps}\partial_x^3\log h_{\eps}\right]\,dx.
\lb{AEnE}
\eea

Given $\eps>0$ equation \rf{ASSM1} is parabolic with respect to $u$. Arguing as in~\citet{BD06}, the initial-boundary value problem system \rf{ASSM1}--\rf{ASSM2} with \rf{ASSM BC}--\rf{ASSM IC} has a unique classical solution at least locally in time. The next proposition establishes a uniform lower bound for $h_{\eps}$ and thus guarantees the global in time solvability for \rf{ASSM1}--\rf{ASSM2}.
\begin{proposition}
For fixed positive $\eps,\sigma,\R,\beta$ and $T>0$ let $(h_{\eps},u_{\eps})$ be the solution of the initial-boundary value problem system \rf{ASSM1}--\rf{ASSM2} with \rf{ASSM BC}--\rf{ASSM IC} in $(0,1)\times (0,T)$. There exists a positive constant $C_3\ge 1$ depending only on $T$, $\alpha$, $\nu$, $\R$, $\beta$, $\sigma$, $h_0$, and $u_0$ such that, for all sufficiently small $\eps>0$, the following terms are bounded by $C_3$ in respective norms
\bea
\sqrt{h_{\eps}},\,\px \sqrt{h_{\eps}},\,h_\eps^{-3/2},\, \px h_{\eps}, \,\sqrt{h_{\eps}}u_{\eps},\,\sqrt{\eps}\partial_x^3 h_{\eps}\in L^\infty(0,T;L^2(0,1)),\lb{AAEs1}\\
\px\left(h_\eps^{-3/2}\right),\,\pxx h_{\eps},\,\sqrt{h_{\eps}}\px u_{\eps},\,u_{\eps},\,\sqrt{\eps}\partial_x^4 h_{\eps},\,\eps\pxx u_{\eps}\in L^2((0,1)\times (0,T)),
\lb{AAEs}
\eea
and 
\be
\frac{1}{C_3}\le h_{\eps}(x,t)\le C_3\,, \quad (x,t)\in (0,1)\times (0,T)\,.
\lb{ABEs}
\ee
\lb{Pr3}
\end{proposition}
\no\textbf{Proof:} Throughout the proof $C$ denotes a positive constant depending on the same parameters as $C_3$ that
  may vary from line to line. Proceeding as in the proof of Proposition~\ref{Pr2}, the mass conservation and the energy equality \rf{AEE} imply the estimate \rf{ABEs} and all bounds in \rf{AAEs1}--\rf{AAEs} except those on $\px \sqrt{h_{\eps}},\,\pxx h_{\eps},\,\px\left(h_\eps^{-3/2}\right)$ and $\sqrt{\eps}\partial_x^4 h_{\eps}$. Using the previously obtained bounds (in particular the lower bound \rf{ABEs} on $h_\eps$) and the interpolation inequalities
$$
\|\pxx h_\eps\|_{L^\infty(0,1)} \le C\, \|\px^3 h_\eps\|_{L^2(0,1)} \ \text{and} \ \|\px h_\eps\|_{L^6(0,1)} \le C\, \|\pxx h_\eps\|_{L^2(0,1)}^{1/3} \, \|\px h_\eps\|_{L^2(0,1)}^{2/3}\,, 
$$  
one can estimate the last term on the right-hand side of the entropy equality \rf{AEnE} as follows:
\bea
&&\left| 4\nu\eps^2\int_0^T\int_0^1\pxx u_\eps\cdot\partial_x^3\log h_\eps\,dx\,dt \right| \nonumber\\
&\le& 4\nu\eps\, \left( \int_0^T \int_0^1 \eps^2 |\pxx u_\eps|^2\,dx\,dt \right)^{1/2}\ \sup_{[0,T]} ||\partial_x^3\log h_\eps||_{L^2(0,1)}\nonumber\\
&\le& C\,\eps\, \sup_{[0,T]} \left\| \frac{\px^3 h_\eps}{h_\eps}-3\frac{\px  h_\eps\pxx h_\eps}{h_\eps^2}+2\frac{|\px h_\eps|^3}
{h_\eps^3}\right\|_{L^2(0,1)}\nonumber\\[2ex]
&\le&C\,\eps\, \sup_{[0,T]}\left\{ \|\px^3 h_\eps\|_{L^2(0,1)}+\|\px h_\eps\|_{L^2(0,1)}\, \|\pxx h_\eps||_{L^\infty(0,1)}
+\|\px h_\eps\|_{L^6(0,1)}^3\right\} \nonumber\\
&\le& C\,\eps\, \sup_{[0,T]}\left\{ \|\px^3 h_\eps\|_{L^2(0,1)} + C\, \|\px^3 h_\eps\|_{L^2(0,1)} + C\, \|\pxx h_\eps\|_{L^2(0,1)} \right\}\nonumber\\
&+&C\,\sqrt{\eps} \, \sup_{[0,T]} \| \sqrt{\eps} \, \px^3 h_\eps\|_{L^2(0,1)} \le C\,\sqrt{\eps}.
\lb{RT}
\eea
Therefore, taking $\eps$ sufficiently small, the bounds in \rf{AAEs1}--\rf{AAEs} still to be proved follow from the entropy inequality \rf{AEnE} exactly as in the proof of Proposition~\ref{Pr2}.$\hfill\blacksquare$

We now show that solutions to the system \rf{ASSM1}--\rf{ASSM2} with boundary and initial conditions \rf{ASSM BC}--\rf{ASSM IC} converge to a solution to \rf{WF1}--\rf{WF2} as $\eps\go 0$.
\begin{theorem}
For any positive $\sigma,\R,\beta$ and initial data $(h_0,u_0)$ satisfying
\rf{IAs}, there exists a global weak solution $(h,u)$ to the system \rf{SSM1}--\rf{SSM2} with boundary conditions \rf{SSM flux BC}--\rf{SSM height BC} and initial conditions $(h_0,u_0)$ in the sense of Definition~\ref{defws}.
\lb{T1}
\end{theorem}
\no\textbf{Proof:}  Take a sequence $\{\eps_n\}_{n\ge 1}\go 0$ and, for each $n\ge 1$, denote the corresponding solution to the approximate system \rf{ASSM1}--\rf{ASSM2}--\rf{ASSM BC}--\rf{ASSM IC} with $\eps=\eps_n$ by $(h_{\eps_n},u_{\eps_n})$. 

We investigate the compactness properties of the sequence $(h_{\eps_n},u_{\eps_n})$ and to this end use the estimates derived in Proposition~\ref{Pr3}. Let $T>0$. First, owing to \rf{AAEs1} and \rf{ABEs}, $(\sqrt{h_{\eps_n}})$ and $(\sqrt{h_{\eps_n}} u_{\eps_n})$ are bounded in $L^\infty((0,1)\times (0,T))$ and $L^\infty(0,T;L^2(0,1))$, respectively. Hence, by \rf{ASSM2}, $(\pt h_{\eps_n})$ is bounded in $L^\infty(0,T;H^{-1}(0,1))$. Since $(h_{\eps_n})$ is bounded in $L^\infty(0,T;H^1(0,1))$ and $L^2(0,T;H^2(0,1))$ by \rf{AAEs1} and \rf{AAEs}, it follows from the compactness of the embedding of $H^1(0,1)$ in $C([0,1])$ and Corollary~4 in~\citet{Si87} that there is $h\in C([0,1]\times [0,T])$ such that, after possibly extracting a subsequence,
\bea
h_{\eps_n}&\go& h\ \ \text{in}\ \ L^2(0,T;W^{1,p}(0,1))\cap C([0,1]\times [0,T])\ \ \text{for}\ \ p\in [1,\infty),\lb{hBd0}\\
\pt h_{\eps_n}&\stackrel{\star}{\rightharpoonup} & \pt h\ \ \text{in}\ \ L^\infty(0,T;H^{-1}(0,1)). \nonumber
\eea
Combining \rf{ABEs} and \rf{hBd0} we additionally obtain that $h$ is positive and
\be
h_{\eps_n}^{-1}\go h^{-1}\ \ \text{in}\ \ C([0,1]\times [0,T]).
\lb{hmBd}
\ee

We next turn to compactness properties of $(u_{\eps_n})$. For that purpose, we write \rf{ASSM1} as
\bea
\R\,\pt(h_{\eps_n}u_{\eps_n})&=&-\R\,\px(h_{\eps_n}u_{\eps_n}^2)+4\px(\nu\,h_{\eps_n}\px u_{\eps_n})
+h_{\eps_n}\px(\sigma\pxx h_{\eps_n}-\Pi(h_{\eps_n}))\nonumber\\
&-&\frac{u_{\eps_n}}{\beta}+\eps_n h_{\eps_n}\partial_x^7 h_{\eps_n}-\eps_n^2\partial_x^4 u_{\eps_n}
\lb{uBd}
\eea
and claim that the right-hand side of \rf{uBd} is bounded in $L^2(0,T;H^{-3}(0,1))$. Indeed, by \rf{AAEs1}--\rf{ABEs}, $(h_{\eps_n} u_{\eps_n}^2)$, $(\sqrt{h_{\eps_n}})$, and $(\sqrt{h_{\eps_n}}\px u_{\eps_n})$ are in $L^\infty(0,T;L^1(0,1))$,  $L^\infty((0,1)\times (0,T))$, and $L^2((0,1)\times (0,T))$, which imply that $(h_{\eps_n} u_{\eps_n}^2)$ and $(\px (h_{\eps_n}\px u_{\eps_n}))$ are bounded in $L^2(0,T;H^{-1}(0,1))$. Next, for any $\psi\in C_0^\infty((0,1)\times (0,T))$ one obtains, using integration by parts and \rf{AAEs1}--\rf{ABEs}, that
\beas
&&\left| \int_0^T\int_0^1\psi h_{\eps_n}\px^3 h_{\eps_n}\,dx\,dt \right| = \left| \int_0^T\int_0^1\pxx h_{\eps_n}\left[h_{\eps_n}\px\psi+\psi\px h_{\eps_n}\right]\,dx\,dt \right|\\
&\le&\int_0^T \|\pxx h_{\eps_n}\|_{L^2(0,1)} \left[ \|h_{\eps_n}\|_{L^\infty(0,1)} \|\px \psi\|_{L^2(0,1)} + \|\psi\|_{L^\infty(0,1)} \|\px h_{\eps_n}\|_{L^2(0,1)} \right]\,dt\\
&\le& C\, \left( \int_0^T \|\psi\|_{H^1(0,1)}^2 \,dt \right)^{1/2}
\eeas
and
\beas
&&\left| \int_0^T\int_0^1\psi\, h_{\eps_n}\px\Pi(h_{\eps_n}) \,dx\,dt \right| = \left| \int_0^T\int_0^1\px\psi\, \Pi_1(h_{\eps_n}) \,dx\,dt \right|\\
&\le& \|\Pi_1(h_{\eps_n})\|_{L^\infty((0,1)\times (0,T))} \left( \int_0^T \|\psi\|_{H^1(0,1)}^2  \,dt\right)^{1/2},
\eeas
where $\Pi_1(h)$ is defined in \rf{Pi1}. Similarly,
\beas
&&\eps_n\, \left| \int_0^T\int_0^1 \psi h_{\eps_n}\px^7h_{\eps_n}\,dx\,dt \right|\\
&=&\eps_n\, \left| \int_0^T\int_0^1 \px^4h_{\eps_n}\, \left( \psi\, \px^3h_{\eps_n} + 3\px\psi\, \pxx h_{\eps_n} + 3\pxx \psi\, \px h_{\eps_n} + h_{\eps_n}\, \px^3\psi \right) \,dx\,dt \right|\\
&\le&\eps_n\, \int_0^T \|\px^4 h_{\eps_n}\|_{L^2(0,1)}\ \left( \|\psi\|_{L^\infty(0,1)} \|\px^3 h_{\eps_n}\|_{L^2(0,1)} +\|h_{\eps_n}\|_{L^\infty(0,1)} \|\px^3 \psi\|_{L^2(0,1)} \right.\\
&& \left. \qquad +\  3 \|\px \psi\|_{L^\infty(0,1)} \|\pxx h_{\eps_n}\|_{L^2(0,1)} + 3 \|\pxx \psi\|_{L^\infty(0,1)} \|\px h_{\eps_n}\|_{L^2(0,1)} \right) \,dt \\
&\le&C\, \sqrt{\eps_n}\, \int_0^T \|\px^4 h_{\eps_n}\|_{L^2(0,1)} \|\psi\|_{H^3(0,1)}\,dt \le C\, \|\psi\|_{L^2(0,T;H^3(0,1))}.
\eeas
Finally, $(u_{\eps_n})$ and $(\eps_n \pxx u_{\eps_n})$ are bounded in $L^2((0,1)\times (0,T))$ and $L^2(0,T;H^{-2}(0,1))$, respectively, by \rf{AAEs}. Collecting the above information completes the proof of the boundedness of the right-hand side of \rf{uBd}, whence
\be
(\pt(h_{\eps_n}u_{\eps_n}))\ \text{ is bounded in }\ L^2(0,T;H^{-3}(0,1)).
\lb{spip}
\ee
Now, owing to \rf{AAEs} and \rf{ABEs}, we realize that $(u_{\eps_n})$ is
bounded in $L^2(0,T;H^1(0,1))$. This fact along with \rf{AAEs1}--\rf{ABEs}
allows us to conclude that $(h_{\eps_n}u_{\eps_n})$ is bounded in
$L^2(0,T;H^1(0,1))$. Combining this with \rf{spip}
and Corollary~4 in~\citet{Si87} ensures that $(h_{\eps_n}u_{\eps_n})$ is compact
in $L^2((0,1)\times (0,T))$. Then, thanks to \rf{hmBd}, there exists $u\in L^2(0,T;H^1(0,1))$ such that 
\be
\px u_{\eps_n}\rightharpoonup\px u\ \ \text{in}\ \ L^2((0,1)\times
(0,T))\ \ \text{and}\ \ u_{\eps_n}\go u \ \ \text{in}\ \ L^2((0,1)\times (0,T)).
\lb{umBd}
\ee
The convergences \rf{hBd0}, \rf{hmBd}, and \rf{umBd} then allow to pass to the limit as $n\go \infty$ in the weak formulation of the approximating systems \rf{ASSM1}--\rf{ASSM2}--\rf{ASSM BC}--\rf{ASSM IC} in order to get that $(h,u)$ satisfies \rf{WF1}--\rf{WF2}.$\hfill\blacksquare$

\begin{remark}[Energy and entropy inequalities]
It follows from the proof of Theorem~\ref{T1} together with the energy and entropy identities \rf{AEE}, \rf{AEnE}, the property \rf{AP} and the estimate \rf{RT} that the weak solutions to \rf{WF1}--\rf{WF2} obtained above satisfy the corresponding energy inequality
\bea
\displaystyle &&\sup_{t\in (0,T]}\int_0^1 \left[ \R\, h\frac{u^2}{2} + U(h) +\sigma\frac{|\px h|^2}{2} \right] \,dx + \int_0^T\int_0^1 \left( 4\nu h|\px u|^2 + \frac{u^2}{\beta} \right) \,dxdt \nonumber\\
&\le& \int_0^1 \left[ \R\, h_0\frac{u_0^2}{2} + U(h_0) + \sigma\frac{|\px   h_0|^2}{2} \right] \,dx
\lb{EI}
\eea
and entropy inequality
\bea
\displaystyle &&\sup_{t\in (0,T]} \int_0^1 \left[ \frac{1}{2}h\, (\R\, u+\px\varphi(h))^2 -\frac{1}{\beta}\,\varphi(h) + \R\,\left( \sigma\frac{|\px h|^2}{2}+U(h) \right) \right] \,dx\nonumber\\
&+&\int_0^T\int_0^1 \left( \R\,\frac{u^2}{\beta} + 4\sigma\nu|\pxx h|^2 + 4\nu\Pi'(h)|\px h|^2 \right) \,dxdt \nonumber\\
&\le&\int_0^1 \left[ \frac{1}{2} h_0\, (\R\, u_0+\px\varphi(h_0))^2 - \frac{1}{\beta}\,\varphi(h_0) + \R\ \left( \sigma\frac{|\px h_0|^2}{2} + U(h_0) \right) \right] \,dx. \phantom{000000}
\lb{EnI}
\eea
Note that the right-hand side of \rf{EI}--\rf{EnI} does not depend on $T$. Consequently, the statements of Lemma~\ref{C1} and Proposition~\ref{Pr2} hold also for the weak solutions to \rf{WF1}--\rf{WF2} constructed in Theorem~\ref{T1}.
\lb{R1}
\end{remark}

\section{Limiting cases.}
\label{S.limit1}

\subsection{Cases $\R=0$ and $\beta=\infty$.}

By Remark~\ref{R1}, the statement of Proposition~\ref{Pr2} is true for the
weak solutions to \rf{WF1}--\rf{WF2} provided by Theorem~\ref{T1}. We may then
investigate the behaviour of these solutions as either $\R\go 0$ or
$\beta\go\infty$. Let us first consider a sequence of positive real numbers
$(\beta_n)$, $\beta_n\go\infty$, and denote the corresponding solutions to
\rf{WF1}--\rf{WF2} with a fixed $\R>0$ and $\beta=\beta_n$ by $(h_{\beta_n},\,u_{\beta_n})$. Though the estimate on $(u_{\beta_n}/\beta_n)$ is useless in that case, one still recovers the estimate of $(u_{\beta_n})$ in $L^2(0,T;H^1_0(0,1))$ as a consequence of \rf{AEs}, \rf{BEs}, and the Poincar\'e inequality.
Arguing as in the proof of Theorem~\ref{T1}, we conclude that, after possibly extracting a subsequence, $(h_{\beta_n},\,u_{\beta_n})$ converges towards a weak solution to the strong-slip model \rf{FFM1}--\rf{FFM2} describing the dynamics of free suspended films. 

If we now consider a sequence of positive real numbers $(\R_n)$, $\R_n\go 0$, and denote the corresponding weak solutions
to \rf{WF1}--\rf{WF2} with $\R=\R_n$ and a fixed $\beta\in (0,\infty]$ by $(h_{\R_n},\,u_{\R_n})$, we may again proceed as in the proof of Theorem~\ref{T1} to show that, after possibly extracting a subsequence, $(h_{\R_n},\,u_{\R_n})$ converges towards a weak solution to the strong-slip model \rf{SSM1}--\rf{SSM2} without inertial terms. There is however an important difference as the bound on $(\pt(\,h_{\R_n}u_{\R_n}))$ in $L^2(0,T;H^{-3}(0,1))$ is no longer available and we only obtain the weaker conclusion
\bes
u_{\R_n}\rightharpoonup u\ \ \text{in}\ \ L^2(0,T;H^1(0,1)).
\ees
Still, owing to the strong convergence of $(h_{\R_n})$, this allows to pass to the limit as $n\to\infty$ in all the terms which depends linearly on $u_{\R_n}$, that is, all the terms involving $u_{\R_n}$ except $\R_n h_{\R_n} u_{\R_n}^2$. But the latter converges to zero as $\R_n\go 0$ since $(h_{\R_n})$ and $(u_{\R_n})$ are bounded in $L^\infty((0,1)\times (0,T))$ and $L^2(0,T;H^1_0(0,1))$, respectively, by \rf{AEs}--\rf{BEs}.

\subsection{Case $\sigma=0$.}

As already pointed out in Remark~\ref{R2}, for weak solutions $(h,u)$ to \rf{WF1}--\rf{WF2}, the estimate on $\px h$ in \rf{AEs} and both bounds in \rf{BEs} depend on $\sigma$ and explode as $\sigma\go 0$.  In order to investigate the limiting behaviour of weak solutions to \rf{WF1}--\rf{WF2} as $\sigma\to 0$, we refine the a priori bounds of Proposition~\ref{Pr2} in the next proposition so as to avoid their dependence on $\sigma$, although paying for it with an exponential growth with respect to time.

\begin{proposition}
For fixed positive $\sigma,\R,\beta,T$  and initial data satisfying \rf{IAs}, let $(h,u)$ be a weak solution to \rf{WF1}--\rf{WF2}. There exists a positive constant $C_4>1$ depending only on $T$, $\alpha$, $\nu$, $\R$, $E(u_0,h_0)$, and $S_\beta(u_0,h_0)$ defined in \rf{Sb} (but not on $\sigma\in (0,1)$) such that the following terms are bounded by $C_4$ in respective norms
\bea
\sqrt{h},\,\px \sqrt{h},\,h^{-3/2},\,\sqrt{h}u&\in& L^\infty(0,T;L^2(0,1)),\lb{volvic}\\
\px(h^{-3/2}),\,\sqrt{\sigma} \pxx h,\,\sqrt{h}\px u,\, u &\in& L^2((0,1)\times (0,T)), \lb{vittel}
\eea
and
\be
C_4^{-1} \le h(x,t)\le C_4\,, \quad (x,t)\in (0,1)\times (0,T).\lb{evian}
\ee
\lb{Pr4}
\end{proposition}

\no\textbf{Proof:} 
First, the $L^\infty(0,T;L^2(0,1))$-estimates on $\sqrt{h}$, $\sqrt{h} u$ and $h^{-3/2}$ and the $L^2((0,1)\times (0,T))$-estimates on $\sqrt{h} \px u$ and $u$ readily follow from the mass conservation and the energy inequality \rf{EI} as before. Next, we actually estimate in a different way the term $4\nu \Pi'(h)|\px h|^2\,dxdt$ in the entropy inequality \rf{EnI}. More precisely, using the estimate
\beas
\Pi'(h)&=&\frac{4\alpha-3h}{h^5}\ge\chi_{(0,\alpha)}(h)\frac{\alpha}{h^5}+\chi_{(\alpha,\infty)}(h)\left(-\frac{3}{h^4}\right)\\
&\ge&\frac{\alpha}{h^5}-\chi_{(\alpha,\infty)}(h)\left(\frac{\alpha}{h^5}+\frac{3}{h^4}\right)\\
&\ge&\frac{\alpha}{h^5}-\chi_{(\alpha,\infty)}(h)\frac{4}{h^4},
\eeas
where $\chi_{A}$ denotes the characteristic function of a set $A\subset\mR$ and recalling that $\varphi(h)=4\nu\log{h}$, one obtains
\be
\Pi'(h)|\px h|^2 \ge \frac{\alpha}{h^5}|\px h|^2-\chi_{(\alpha,\infty)}(h)\frac{h|\px\varphi(h)|^2}{4\nu^2\alpha^3}\ge\frac{\alpha}{h^5}|\px h|^2-\frac{h|\px\varphi(h)|^2}{4\nu^2\alpha^3}.
\lb{PiEs}
\ee
Next, let us consider as in the proof of Proposition \ref{Pr2} the function $S_\beta(u,h)$  defined in \rf{Sb}. Thanks to the mass conservation \rf{conservation of mass}, the entropy inequality \rf{EnI} also reads
\bea
S_\beta(u(t),h(t)) &+& \int_0^t\int_0^1\left(\R\frac{u^2}{\beta}+4\sigma\nu|\pxx h|^2\right)\,dx\,ds \nonumber\\
&\le& S_\beta(u_0,h_0)-4\nu\int_0^t \int_0^1\Pi'(h)|\px h|^2\,dx\,ds. \lb{SEqs}
\eea
Recalling that
\be
S_\beta(u,h) \ge \frac{1}{2}\int_0^1h(\R\, u+\px\varphi(h))^2\,dx-\frac{\R}{6\alpha^2}
\lb{SIns}
\ee
by \rf{contrex}, it follows from \rf{PiEs}, \rf{SIns}, and Lemma~\ref{C1} that
\beas
& & -4\nu\int_0^1\Pi'(h)|\px h|^2\,dx \le -4\alpha\nu\int_0^1\frac{|\px h|^2}{h^5}\,dx+\frac{1}{\nu\alpha^3}\int_0^1h|\px\varphi(h)|^2\,dx\\
&\le&-4\alpha\nu\int_0^1\frac{|\px h|^2}{h^5}\,dx+\frac{4}{\nu\alpha^3}\,\left(S_\beta(u,h)+\R\,E(u,h)+\frac{\R}{3\alpha^2} \right) \\
&\le&-4\alpha\nu\int_0^1\frac{|\px h|^2}{h^5}\,dx+\frac{4}{\nu\alpha^3}\,\left( S_\beta(u,h)+\R\,E(u_0,h_0)+\frac{\R}{3\alpha^2} \right),
\eeas
where we used the energy inequality \rf{EI} in the last estimate. Inserting this in \rf{SEqs} one ends up with
\beas
&& S_\beta(u(t),h(t))+\int_0^t\int_0^1\left(\R\frac{u^2}{\beta}+4\sigma\nu|\pxx h|^2 +4\alpha\nu\frac{|\px h|^2}{h^5} \right)\,dx\,ds \\
&\le& S_\beta(u_0,h_0) + \frac{4\R}{\nu\alpha^3}\, \left(E(u_0,h_0) + \frac{1}{3\alpha^2} \right)\, t+\frac{4}{\nu\alpha^3} \int_0^t S_\beta(u(s),h(s)) \,ds.
\eeas
Applying Gronwall's inequality gives
\bea
&&S_\beta(u(t),h(t))+\int_0^t\int_0^1\left( \R\frac{u^2}{\beta}+4\sigma\nu|\pxx h|^2\right)\,dx\,dt\nonumber\\ 
&\le&\left(S_\beta(u_0,h_0)+\R\left[ E(u_0,h_0)+\frac{1}{3\alpha^2}\right]\right)\exp\left(\frac{4t}{\nu\alpha^3}\right).
\lb{exB}
\eea
Combining \rf{CoEs}, \rf{SIns}, and \rf{exB} implies the estimates on $\px\sqrt{h}$ in $L^\infty(0,T;L^2(0,1))$ and $\px(h^{-3/2})$, $\sqrt{\sigma}\pxx h$ in $L^2((0,1)\times(0,T))$ and completes the proof of \rf{volvic}--\rf{vittel}. To check \rf{evian}, we first notice that the bound on $\sqrt{h}$ in $L^\infty(0,T;H^1(0,1))$ from \rf{volvic} and the embedding of $H^1(0,1)$ in $L^\infty(0,1)$ guarantee the upper bound in \rf{evian}. Next, integrating the following identity
\bes
\frac{1}{h(x)}=\frac{1}{h(y)}-\int_y^x\frac{\partial_z h(z)}{h^2(z)}\,dz\,, \quad x\in (0,1)\,,
\ees
 with respect to $y$ and using the Cauchy-Schwarz inequality give
\beas
\frac{1}{h(x)} &\le& \int_0^1\frac{dy}{h(y)} + \left(\int_0^1\frac{|\partial_z h(z)|^2}{h(z)}\,dz\right)^{1/2}\, 
\left(\int_0^1\frac{dz}{h^3(z)} \right)^{1/2}\\
&\le& \|h^{-3/2}\|_{L^2(0,1)}^{2/3} + 4\,\|\px\sqrt{h}\|_{L^2(0,1)}\,\|h^{-3/2}\|_{L^2(0,1)}.
\eeas
The lower bound in \rf{evian} then follows from the above inequality by \rf{volvic}.$\hfill\blacksquare$

Using the a priori bounds from Proposition~\ref{Pr4} one can show existence of weak solutions to \rf{WF1}--\rf{WF2} in the case $\sigma=0$. 
\begin{theorem}
Let $\sigma=0$. For any positive $\R,\beta$ and initial data $(h_0,u_0)$ satisfying \rf{IAs} there exists a global weak solution $(h,u)$ to \rf{SSM1}--\rf{SSM2} satisfying \rf{WF1}--\rf{WF2}.
\lb{T3}
\end{theorem}
\no\textbf{Proof:}
Let us take a sequence $\{\sigma_n\}\go 0$ and, for each $n\ge 1$, denote the corresponding weak solution to \rf{WF1}--\rf{WF2}  with $\sigma=\sigma_n$ by $(h_{\sigma_n},u_{\sigma_n})$. Owing to Proposition~\ref{Pr4} we may proceed as in the proof of Theorem~\ref{T1} to show the existence of functions $h$ and $u$ satisfying
\beas
h_{\sigma_n}&\go& h\ \ \text{and} \ \ h_{\sigma_n}^{-1} \go  h^{-1} \ \text{in}\ \ C([0,1]\times [0,T]),\\
\pt h_{\sigma_n}&\stackrel{*}{\rightharpoonup}& \pt h\ \ \text{in}\ \ L^\infty(0,T;H^{-1}(0,1))
\eeas 
and
\bes
\px u_{\sigma_n}\rightharpoonup\px u\ \ \text{in}\ \ L^2((0,1)\times (0,T))\ \ \text{and}\ \ u_{\sigma_n}\go u \ \ \text{in}\ \ L^2((0,1)\times (0,T)).
\ees
In addition, 
$$
\sigma_n \pxx h_{\sigma_n} \go 0 \ \ \text{in}\ \ L^2((0,1)\times (0,T))
$$
by \rf{vittel}, from which we deduce that 
$$
\sigma_n \px h_{\sigma_n} \pxx h_{\sigma_n} \rightharpoonup 0 \ \ \text{and}\ \ \sigma_n h_{\sigma_n} \pxx h_{\sigma_n} \rightharpoonup 0 \ \ \text{in}\ \ L^1((0,1)\times (0,T))
$$
with the help of \rf{volvic}. This allows us to pass to the limit as $n\go \infty$ in all integrals in
\rf{WF1}--\rf{WF2} and conclude that $(h,u)$ satisfies \rf{WF1}--\rf{WF2} with $\sigma=0$.$\hfill\blacksquare$ 
  
\section{Convergence to solutions of intermediate-slip \\ equation.}
\label{S.limit2}

In the introduction we mentioned that performing the change of variables \rf{CoV} and passing formally to the limit as $\beta\go 0$ in \rf{SSM1}--\rf{SSM2} gives\rf{ISM1}--\rf{ISM2} and consequently that solutions of the strong-slip equation are expected to converge to that of \rf{ISM}. The next theorem makes this formal procedure rigorous. More precisely, we show that scaled weak solutions to \rf{SSM1}--\rf{SSM2} converge to that of  \rf{ISM} when $\sigma> 0$. In this section, we denote the unscaled variables in \rf{SSM1}--\rf{SSM2} with bars.
\begin{theorem}
For fixed positive $\R,\,\sigma,$ let $(\beta_n)$ be a sequence of positive
real numbers, $\beta_n\go 0$, and, for $n\ge 1$, denote a global weak solution
to \rf{WF1}--\rf{WF2} with initial data $(h_0,u_0)$ satisfying \rf{IAs} and  $\beta=\beta_n$ by $(\bar{h}_n,\bar{u}_n)$. Using \rf{CoV}, we define
\be
h_n(x,t):=\bar{h}_n\left(x,\frac{t}{\beta}\right) \ \ \text{and}\ \ u_n(x,t):=\frac{1}{\beta}\,\bar{u}_n\left(x,\frac{t}{\beta}\right)\,, \quad (x,t)\in (0,1)\times (0,\infty).
\lb{CoF}
\ee
Then there exist a positive function $h$ and a subsequence of $(h_n,u_n)$ (not relabeled) such that, for any $T>0$,
\beas
&&h_n\go h\ \ \text{in}\ \ L^2(0,T;H^1(0,1))\cap C([0,1]\times [0,T]),\nonumber\\
&&u_n\rightharpoonup u:=h\px(\sigma\pxx h-\Pi(h))\ \ \text{in}\ \ L^2((0,1)\times (0,T)), \lb{CSI}
\eeas
and $h$ is a global weak solution to \rf{ISM} considered with initial condition $h_0$ satisfying \rf{IAs} and boundary conditions \rf{SSM flux BC}--\rf{SSM height BC}. \lb{T2}
\end{theorem}
\no\textbf{Proof:} Owing to \rf{WF1}--\rf{WF2}, for each $n\ge 1$, the functions $(h_n,u_n)$ satisfy the following scaled weak formulation
\sbea
\lb{WF1s}
&&\int_0^\infty\int_0^1h_n\pt\psi\,dx\,dt + \int_0^1 h_0 \psi(.,0)\,dx=-\int_0^\infty\int_0^1h_nu_n\px\psi\,dxdt,\\[2ex]
&&\beta_n^2\left(\R\,\int_0^\infty\int_0^1h_nu_n\pt\phi\,dx\,dt +\R\,\int_0^1h_0u_0\phi(.,0)\,dx + \R\,\int_0^\infty\int_0^1 h_nu_n^2\px\phi \,dx\,dt\right)\nonumber\\
&& -\,4\beta_n\nu\int_0^\infty\int_0^1h_n\px u_n\px \phi\,dx\,dt - \sigma\int_0^\infty\int_0^1\px h_n\pxx h_n\phi\,dx\,dt
\nonumber\\
&& -\, \int_0^\infty\int_0^1 \left( \sigma\,h_n\pxx h_n -\Pi_1(h_n) \right) \px\phi\,dx\,dt - \int_0^\infty\int_0^1u_n\phi\,dx\,dt=0
\lb{WF2s}
\seea
for all $\psi\in C_0^\infty([0,1]\times [0,\infty))$ and $\phi\in C_0^\infty((0,1)\times [0,\infty))$. In addition, for any $T>0$, $(h_n,u_n)$ satisfies the scaled energy inequality
\beas
& & \sup_{t\in[0,T]}\int_0^1 \left[ \beta_n^2\R\,h_n\frac{u_n^2}{2}+U(h_n) +\sigma\frac{|\px h_n|^2}{2} \right] \,dx + \int_0^T\int_0^1 \left( 4\nu\beta_n\, h_n|\px u_n|^2 + u_n^2 \right) \,dx\,dt \\
& &\le \int_0^1 \left[ \beta_n^2\R\, h_0\frac{u_0^2}{2}+U(h_0)+\sigma\frac{|\px h_0|^2}{2} \right]\,dx
\eeas
and the scaled entropy inequality
\beas
& &\sup_{t\in[0,T]}\int_0^1 \left[ \frac{\beta_n}{2}h_n(\beta_n\R\,u_n+\px\varphi(h_n))^2 +\beta_n\R\left(\sigma\frac{|\px h_n|^2}{2}+U(h_n)\right)- \varphi(h) \right] \,dx \\
& & +\, \int_0^T\int_0^1 \left( \beta_n\R\, u_n^2+4\sigma\nu|\pxx h_n|^2+4\nu\Pi'(h_n)|\px h_n|^2 \right)\,dx\,dt\\
& &\le \beta_n\int_0^1 \left[ \frac{1}{2} h_0(\beta_n\R\, u_0+\px\varphi(h_0))^2+\R \left( \sigma\frac{|\px h_0|^2}{2}+U(h_0) \right) \right] \,dx - \int_0^1 \varphi(h_0)\,dx.
\eeas
Using the fact that $\sigma>0$ and proceeding as in the proof of Proposition~\ref{Pr2}, it follows from the mass conservation and the scaled energy inequality that 
\begin{eqnarray}
& & (\sqrt{h_n}),\, (h_n^{-3/2}),\,(\px h_n)\ \ \text{are bounded in}\ \ L^\infty(0,T;L^2(0,1)), \lb{IAEs1}\\
& & (u_n),\, (\sqrt{\beta_n} \sqrt{h_n} \px u_n)\ \ \text{are bounded in}\ \ L^2((0,1)\times (0,T)), \lb{IAEs2}
\end{eqnarray}
and there exists $C_5>1$ such that 
\bes
C_5^{-1}\le h_n(x,t)\le C_5\,, \quad (x,t)\in (0,1)\times (0,T).
\lb{IBEs}
\ees
We next use the scaled entropy inequality as at the end of the proof of Proposition~\ref{Pr2} to establish that 
\bea
(\px(h_n^{-3/2})), \, (\pxx h_n), \ \ \text{are bounded in}\ \ L^2((0,1)\times (0,T)).
\lb{IAEs}
\eea
Proceeding as in the proof of Theorems~\ref{T1} and~\ref{T3} one obtains that there exist functions $h$ and $u$ such that
\beas
&&h_n\go h\ \ \text{in}\ \ L^2(0,T;H^{1}(0,1))\cap C([0,1]\times[0,T]),\nonumber\\
&&h_n^{-1}\go h^{-1}\ \ \text{in}\ \ C([0,1]\times[0,T]),\nonumber\\
&&u_n\rightharpoonup u\ \ \text{in}\ \ L^2((0,1)\times (0,T)), \nonumber \\
&& \beta_n h_n \px u_n \go 0 \ \ \text{in}\ \ L^2((0,1)\times (0,T)). \nonumber
\eeas
It is then straightforward to pass to the limit as $\beta_n\to 0$ in \rf{WF1s}--\rf{WF2s} to obtain
\sbea
\lb{IWF1}
& &\int_0^\infty\int_0^1h\pt\psi\,dx\,dt + \int_0^1 h_0 \psi(.,0)\, dx = -\int_0^\infty\int_0^1hu\px\psi\,dx\,dt,\\
& & -\sigma\int_0^\infty\int_0^1\pxx h\px (h\phi)\,dxdt + \int_0^\infty \int_0^1 \Pi_1(h) \px\phi\,dxdt = \int_0^\infty \int_0^1u\phi\,dx\,dt. \phantom{000}
\lb{IWF2}
\seea
Using an approximation argument, we may actually take $\phi=h\px\psi$ in \rf{IWF2} with $\psi\in C_0^\infty((0,1)\times [0,\infty))$ and obtain
\beas
\int_0^\infty \int_0^1h\pt\psi\,dx\,dt + \int_0^1 h_0 \psi(.,0)\, dx & = & \sigma\int_0^\infty\int_0^1\pxx h\px (h^2\px \psi)\,dx\,dt \\
& & +\ \int_0^\infty \int_0^1h^2\px(\Pi(h))\psi\,dx\,dt,
\eeas
which coincides with the weak formulation for \rf{ISM} introduced in~\citet{BF90} and~\citet{BP98}. $\hfill\blacksquare$

\medskip

\begin{remark}
Owing to the positivity of the weak solution $h$ to \rf{ISM} constructed in
Theorem~\ref{T2}, equation \rf{ISM} is then uniformly parabolic and it is
likely that classical parabolic regularity results ensure that $h$ is actually
a classical solution to \rf{ISM}. As it is shown in~\citet{BGW01} that there
is a unique positive classical solution to \rf{ISM}, this implies the
convergence of the whole sequence $(h_n)$ towards it.
\end{remark}

\section{Discussion.}

We would like to note some open questions  that  remain to be solved for the strong-slip system \rf{SSM1}--\rf{SSM2} as well as the corresponding shallow-water equations \rf{GKM1}--\rf{GKM2}. First, it would be of interest to figure out whether the weak solutions constructed in the previous sections are more regular. Partial answers have already been given for particular cases of \rf{GKM1}--\rf{GKM2} with $\sigma=0$, see~\citet{So76} and~\citet{MV08} where the existence of strong solutions is established. Furthermore, in~\citet{VLW01}  the  two-dimensional  analogue of \rf{FFM1}--\rf{FFM2} describing dewetting of three-dimensional free suspended films is derived and reads
\begin{subequations}
\begin{align}
\pt(hu)+\dv(hu\otimes u)&=\dv(h(\nabla u+(\nabla u)^T))+2\nabla(h\dv u)+\sigma h \nabla\Delta h-\nabla P(h)
\label{2FFM1}\\
\pt h &= -\,\dv\left(h u\right),
\label{2FFM2}
\end{align}
\end{subequations} 
which corresponds to the choice $\mu(h)=h$ and $\lambda(h)=2h$ in \rf{GKM1}--\rf{GKM2}. As mentioned in the introduction the BD entropy equality has been so far only established for \rf{GKM1}--\rf{GKM2} if the relation \rf{lmR} for the viscosities is satisfied. Unfortunately this relation is not fulfilled for \rf{2FFM1}--\rf{2FFM2}. Analogously a two-dimensional version of \rf{SSM1}--\rf{SSM2} can be constructed but will present the same difficulty, so that the analysis presented herein is not likely to extend in a straightforward way to the two-dimensional model.

\section{Acknowledgements} 

GK acknowledges the support from the Weierstrass Institute and the postdoctoral scholarship at the
Max-Planck-Institute for Mathematics in the Natural Sciences, Leipzig. 
PhL and BN were partially supported by the CNRS/RSE project JP090230 and the EPSRC Science and Innovation award to the
Oxford Centre for Nonlinear PDE (EP/E035027/1).

\bibliographystyle{unsrtnat}
\bibliography{bibliography}

\begin{thebibliography}{20}
\providecommand{\natexlab}[1]{#1}
\providecommand{\url}[1]{\texttt{#1}}
\expandafter\ifx\csname urlstyle\endcsname\relax
  \providecommand{\doi}[1]{doi: #1}\else
  \providecommand{\doi}{doi: \begingroup \urlstyle{rm}\Url}\fi

\bibitem[Oron et~al.(1997)Oron, Davis, and Bankoff]{ODB97}
A.~Oron, S.~H. Davis, and S.~G. Bankoff.
\newblock Long-scale evolution of thin liquid films.
\newblock \emph{Rev. Mod. Phys.}, 69\penalty0 (3):\penalty0 931--980, 1997.

\bibitem[M\"unch et~al.(2006)M\"unch, Wagner, and Witelski]{MWW06}
A.~M\"unch, B.~Wagner, and T.~P. Witelski.
\newblock Lubrication models with small to large slip lengths.
\newblock \emph{J. Engr. Math.}, 53:\penalty0 359--383, 2006.

\bibitem[de~Gennes(1985)]{gennes85}
P.G. de~Gennes.
\newblock Wetting: Statics and dynamics.
\newblock \emph{Rev. Mod. Phys.}, 57:\penalty0 827--863, 1985.

\bibitem[Fetzer et~al.(2007)Fetzer, M{\"u}nch, Wagner, Rauscher, and
  Jacobs.]{FMWJ07}
R.~Fetzer, A.~M{\"u}nch, B.~Wagner, M.~Rauscher, and K.~Jacobs.
\newblock Quantifying hydrodynamic slip: A comprehensive analysis of dewetting
  profiles.
\newblock \emph{Langmuir}, 23:\penalty0 10559--10566, 2007.

\bibitem[Fetzer et~al.(2005)Fetzer, Jacobs, M\"unch, Wagner, and
  Witelski.]{FJMWW05}
R.~Fetzer, K.~Jacobs, A.~M\"unch, B.~Wagner, and T.~P. Witelski.
\newblock New slip regimes and the shape of dewetting thin liquid films.
\newblock \emph{Phys. Rev. Lett.}, 95:\penalty0 127801, 2005.

\bibitem[M\"unch and Wagner(2005)]{MW05}
A.~M\"unch and B.~Wagner.
\newblock Contact-line instability of dewetting thin films.
\newblock \emph{Physica D}, 209:\penalty0 178--190, 2005.

\bibitem[Kargupta et~al.(2004)Kargupta, Sharma, and Khanna]{KSK04}
K.~Kargupta, A.~Sharma, and R.~Khanna.
\newblock Instability, dynamics and morphology of thin slipping films.
\newblock \emph{Langmuir}, 20:\penalty0 244--253, 2004.

\bibitem[Brenner and Gueyffier(1999)]{BG99}
M.~P. Brenner and D.~Gueyffier.
\newblock On the bursting of viscous films.
\newblock \emph{Phys. Fluids}, 11\penalty0 (3):\penalty0 737--739, 1999.

\bibitem[Bresch et~al.(2003)Bresch, Desjardins, and Lin]{BDL03}
D.~Bresch, B.~Desjardins, and C.-K. Lin.
\newblock On some compressible fluid models: Korteweg, lubrication and shallow
  water systems.
\newblock \emph{Comm. Partial Differential Equations}, 28\penalty0
  (3-4):\penalty0 1009--1037, 2003.

\bibitem[Bresch(2009)]{Br09}
D.~Bresch.
\newblock Shallow-water equations and related topics.
\newblock In C.~M. Dafermos and M.~Pokorny, editors, \emph{Handbook of
  differential equations. Evolutionary equations}, volume~5, pages 1--105.
  Elsevier, 2009.

\bibitem[Solonnikov(1976)]{So76}
V.~A. Solonnikov.
\newblock The solvability of the initial-boundary value problem for the
  equations of motion of a viscous compressible fluid.
\newblock \emph{Zap. Naucn. Sem. Leningrad. Otdel. Mat. Inst. Steklov. (LOMI)},
  56\penalty0 (197):\penalty0 128--142, 1976.

\bibitem[Bresch and Desjardins(2003)]{BD03}
D.~Bresch and B.~Desjardins.
\newblock {Existence of global weak solutions for a 2D viscous shallow water
  model and convergence to the quasigeostrophic model}.
\newblock \emph{Comm. Math. Phys.}, 238\penalty0 (1-2):\penalty0 211--223,
  2003.

\bibitem[Mellet and Vasseur(2008)]{MV08}
A.~Mellet and A.~Vasseur.
\newblock {Existence and uniqueness of global strong solutions for
  one-dimensional compressible Navier-Stokes equations}.
\newblock \emph{SIAM J. Math. Anal.}, 39\penalty0 (4):\penalty0 1344--1365,
  2008.

\bibitem[Bresch and Desjardins(2006)]{BD06}
D.~Bresch and B.~Desjardins.
\newblock {On the construction of approximate solutions for the 2D viscous
  shallow water model and for compressible Navier-Stokes models}.
\newblock \emph{J. Math. Pures Appl.}, 86:\penalty0 362--368, 2006.

\bibitem[Bernis and Friedman(1990)]{BF90}
F.~Bernis and A.~Friedman.
\newblock Higher order nonlinear degenerate parabolic equations.
\newblock \emph{J. Differential Equations}, 83:\penalty0 179--206, 1990.

\bibitem[Kitavtsev and Wagner(2010)]{KW09}
G.~Kitavtsev and B.~Wagner.
\newblock Coarsening dynamics of slipping droplets.
\newblock \emph{J. Engr. Math.}, 66:\penalty0 271--292, 2010.

\bibitem[Bertozzi et~al.(2001)Bertozzi, Gr\"un, and Witelski]{BGW01}
A.~L. Bertozzi, G.~Gr\"un, and T.~P. Witelski.
\newblock Dewetting films: bifurcations and concentrations.
\newblock \emph{Nonlinearity}, 14:\penalty0 1569--1592, 2001.

\bibitem[Simon(1987)]{Si87}
J.~Simon.
\newblock {Compact sets in the space $L^p(0,T;B)$}.
\newblock \emph{Ann. Mat. Pura Appl.}, 146\penalty0 (4):\penalty0 65–--96,
  1987.

\bibitem[Bertozzi and Pugh(1998)]{BP98}
A.~L. Bertozzi and M.~C. Pugh.
\newblock Long-wave instabilities and saturation in thin film equations.
\newblock \emph{Comm. Pure Appl. Math.}, 51:\penalty0 625--661, 1998.

\bibitem[Vaynblat et~al.(2001)Vaynblat, Lister, and Witelski]{VLW01}
D.~Vaynblat, J.~R. Lister, and T.~P. Witelski.
\newblock {Rupture of thin viscous films by van der Waals forces: Evolution and
  self-similarity}.
\newblock \emph{Phys. Fluids}, 13\penalty0 (5):\penalty0 1130--1140, 2001.

\end{thebibliography}
\clearpage
\addtocounter{tocdepth}{2}


\end{document}